\documentclass[a4paper,11pt,leqno]{article}
\pdfoutput=1
\usepackage{euscript}
\usepackage{amsmath,amsfonts,amsthm,amssymb}
\usepackage[left=3cm, right=3cm, top=3cm, bottom=3.5cm]{geometry}
\usepackage{graphicx}
\usepackage{mathrsfs}
\usepackage[final]{pdfpages}
\usepackage{enumerate}

\numberwithin{equation}{section}

\theoremstyle{plain}
\newtheorem{theorem}{Theorem}
\newtheorem{lemma}{Lemma}

\theoremstyle{definition}

\newcommand{\N}{\mathbb{N}}
\newcommand{\Z}{\mathbb{Z}}
\newcommand{\R}{\mathbb{R}}
\newcommand{\norm}[1]{\left\Vert#1\right\Vert}

\def\cA{\EuScript{A}}
\def\cE{\EuScript{E}}
\def\cR{\EuScript{R}}
\def\cP{\EuScript{P}}
\def\cF{\EuScript{F}}
\def\cH{\EuScript{H}}
\def\cC{\EuScript{C}}
\def\cD{\EuScript{D}}
\def\cS{\EuScript{S}}
\def\cT{\EuScript{T}}

\def\cl{\mathop{\mathrm{cl}}\nolimits}
\def\id{\mathop{\mathrm{id}}\nolimits}
\def\spec{\mathop{\mathrm{spec}}\nolimits}
\def\span{\mathop{\mathrm{span}}\nolimits}
\def\Re{\mathop{\mathrm{Re}}\nolimits}

\def\codim{\mathop{\mathrm{codim}}\nolimits}
\def\dim{\mathop{\mathrm{dim}}\nolimits}

\def\otto{%
    \leavevmode\vbox{\offinterlineskip
    \halign{\tabskip=0pt\hfil##\hfil\cr
            \vrule height0.4pt depth0pt width0.6em\cr
            \vrule height 1.6ex \cr
            \noalign{\vskip  -1.3ex}
            $\cap$\cr}}} 

\title{\large Classification of connection graphs of global attractors \\ for $S^1$-equivariant 
parabolic equations}

\author
        {Carlos Rocha\\
        Centro de An\'alise Matem\'atica, Geometria e Sistemas Din\^amicos,\\
				Departamento de Matem\'atica, Instituto Superior T\'ecnico, Universidade de Lisboa\\
        Av. Rovisco Pais 1, 1049--001 Lisbon, PORTUGAL\\
        {\tt crocha@tecnico.ulisboa.pt}\\
        {\tt http://camgsd.tecnico.ulisboa.pt}\\
        {} }

\date{version of \today}

\begin{document}

\maketitle

\setlength{\parindent}{0cm}

\begin{abstract}
We consider the characterization of global attractors $\cA_f$ for semiflows generated by scalar 
one-dimensional semilinear parabolic equations of the form $u_t = u_{xx} + f(u,u_x)$, defined 
on the circle $x\in S^1$, for a class of reversible nonlinearities. 
Given two reversible nonlinearities, $f_0$ and $f_1$, with the same lap signature, we prove the 
existence of a reversible homotopy $f_\tau, 0\le\tau\le 1$, which preserves all heteroclinic 
connections. Consequently, we obtain a classification of the connection graphs of global attractors 
in the class of reversible nonlinearities. 
We also describe bifurcation diagrams which reduce a global attractor $\cA_1$ to the trivial global 
attractor $\cA_0=\{0\}$.

\end{abstract}

\section{Introduction}

Consider the scalar semilinear parabolic equation 
\begin{equation} \label{eq101} 
u_t = u_{xx} + f(u,u_x) \quad , \quad x\in S^1=\R/2\pi\Z \ ,
\end{equation}
where the nonlinearity $f:\R^2\rightarrow \R$ is $C^2$ and dissipative. 
Sufficient dissipative conditions are boundedness, $|f(\cdot,\cdot)|\le C_0$, and a sign condition 
$uf(u,0)<0$ for all large $|u|>K_0$. For less restrictive conditions see \cite{ama85,mana97}.
Then, \eqref{eq101} generates a semiflow in the Sobolev space $X=H^s(S^1), s>\frac{3}{2}$, which 
possesses a nonempty compact {\em global attractor} $\cA\subset X$. For details see \cite{firowo04} 
and also \cite{anfi88, mat88,mana97} for initial references. As general references, see 
\cite{hen81, paz83} for semiflows, and \cite{hal88, bavi92, hamaol02} for global attractors. 

Stationary solutions of \eqref{eq101} satisfy the equation
\begin{equation} \label{eq102} 
0 = v_{xx} + f(v,v_x) \quad , \quad x\in S^1 \ ,
\end{equation}
and are either {\em homogeneous equilibria} $v(t,x)=e$, where $f(e,0)=0$, or {\em nonhomogeneous 
stationary waves}, i.e. $2\pi$-periodic solutions of \eqref{eq102}. 
In general, \eqref{eq101} also features time periodic solutions. These are {\em rotating waves} 
$u(t,x)=v(x-ct)$, rotating around the circle $S^1$ with constant speed $c\ne 0$. Moreover, the 
nonconstant wave shapes $v$ correspond to the $2\pi$-periodic solutions of the ODE
\begin{equation} \label{eq103} 
0 = v_{xx} + f(v,v_x) + cv_x \ .
\end{equation}

Our first assumption is 

\bigskip

{\bf (H)}: all equilibria and periodic orbits of \eqref{eq101} are {\em hyperbolic}. 

\bigskip

Hyperbolicity here concerns the sets of homogeneous equilibria and nonhomogeneous time periodic 
solutions (rotating waves) of the semilinear parabolic PDE \eqref{eq101}. 
These are, respectively, the sets $\cE$ and $\cP$ which are finite due to hyperbolicty (H) and 
the compacteness of $\cA$, \cite{firo00,firowo04}:
\begin{equation} \label{eq104}
\begin{aligned}
\cE &= \{ e_j : f(e_j,0)=0, 1\le j\le n \} \ , \\
\cP &= \{ v_j(x) : v_j(0)=v_j(2\pi), p_j(0)=p_j(2\pi), \ x\in[0,2\pi], \ 1\le j\le q \} \ , 
\end{aligned}
\end{equation} 
where $p=v_x$, $v_j$ denotes the rotating wave profiles, i.e. solutions of \eqref{eq103}, and 
$q$ denotes the number of nonhomogeneous periodic solutions properly $x$-shifted so that 
$v_j(0)=\min_{x\in S^1}v_j(x)$. For simplicity we include in $\cP$ also the nonhomogeneous stationary 
waves.
We point out that the critical solutions may have very large Morse indices $i(e_j)=\dim W^u(e_j)$ 
and $i(v_j)=\dim W^u(v_j)-1$, see \cite{firowo04}. Here $W^u(\cdot)$ denotes the unstable 
manifold of an equilibrium or periodic orbit, \cite{hipush77,hamaol02,hen81}.
In contrast, the equilibria of the related ODE \eqref{eq102}, ${\bf e_j}=(e_j,0)$ and 
${\bf v_j}=(v_j,p_j)$, has only alternating saddles and centers, all in one-to-one correspondence 
with the equilibria of \eqref{eq101}. For example, the first saddle corresponds to the first stable 
equilibrium ${\bf e}_1$ and the last saddle corresponds to the last stable equilibrium ${\bf e}_n$.

Notice that a nonhomogeneous stationary solution $v=v(\cdot)$ is not isolated, since all its shifted 
copies $v(\cdot+\vartheta), \vartheta\in S^1$, are also $2\pi$-periodic solutions of \eqref{eq102}. 
Hence, a nonhomogeneous stationary solution is not hyperbolic in stricto sensu. 
To overcome this inconvenience, hyperbolicity here is understood as {\em normal hyperbolicity}, 
\cite{firowo12a, roc12}, much like in the case of a periodic orbit. 

\bigskip

The global attractor $\cA$ is called a {\em Sturm attractor} (see \cite{firo96}) due to the decay 
property \eqref{eq205} of the {\em zero number}, \cite{mana97}. If all the critical elements in 
$\cE\cup\cP$ are strictly hyperbolic, i.e. $\cP$ does not contain nonhomogeneous stationary 
solutions, the Sturm attractor has the {\em Morse-Smale} property, \cite{czro08, jora10}. 
In this case all the {\em heteroclinic orbit connections} between equilibria and periodic orbits 
of $\cA$, which all together compose the global attractor, are determined by a {\em Sturm permutation}, 
\cite{firowo12b,roc12}. See also \cite{furo91} for the initial results on Sturm permutations. 
Unfortunately, for the reversible class of nonlinearities $f(u,v)=f(u,-v)$ included in the present 
$S^1$-equivariant case, we are missing a proper Morse-Smale Theorem. Therefore, our results must be 
restricted to connection equivalence of global attractors postponing orbit equivalence to further 
research studies. 
This observation also applies to the statement of Theorem 3 in \cite{roc12} which should address 
connection equivalence of global attractors instead of orbit equivalence.

\bigskip

The phase portrait of \eqref{eq102} for a reversible nonlinearity $f(u,v)=f(u,-v)$ (even in $v$) 
has all rotating waves frozen, i.e. they are nonhomogeneous equilibria. In addition, \eqref{eq102} 
is {\em integrable} in the phase plane region corresponding to the {\em cyclicity set} $\cC$, i.e. 
the region of all spatially periodic orbits, see \cite{firowo12a}. Hence, for the reversible 
nonlinearity $f(u,u_x)$ we obtain a {\em period map} $T:\cD\subset\R\rightarrow\R_+$ where the 
domain $\cD$ corresponds to (the $u$-values) of the cyclicity set $\cC$.
 
The period map $T=T(u_0)$, also called time map, is essential for the characterization of all 
$2\pi$-periodic solutions of autonomous planar Hamiltonian equations, see \cite{firowo12a,roc07}. 
In our setting, the period map is used for the characterization of all the equilibria of \eqref{eq101}, 
in particular the frozen rotating waves.

Our objective is the classification of the connection graphs of all Sturm global attractors of flows 
generated by \eqref{eq101}. 
This is achieved using the {\em lap signature class} introduced in \cite{roc12, firowo12a}. 
The lap signature of a period map $T=T(u_0)$ consists of the set of {\em period lap numbers} of the 
$2\pi$-periodic solutions of \eqref{eq102} endowed with a total order derived from the nesting of the 
periodic orbits in the phase space $(v,v_x)$, and their {\em regular parenthesis structure} (called 
regular bracket structure in \cite{lazv04,roc12}). We notice that the period lap number $\ell(v)$ of 
$v\in\cP$ is half of the zero number: $z(v)=2\ell(v)$, see \cite{firowo12a}. The {\em lap signature 
class} of a period map $T=T(u_0)$ is the set of period maps with the same lap signature of $T=T(u_0)$.

\bigskip

Let $\cR$ denote the space of reversible nonlinearities. Our main result asserts the following:

\bigskip

\begin{theorem} \label{th1}
Let $f=f(u,u_x)\in\cR$ and $g=g(u,u_x)\in\cR$ denote two reversible nonlinearities with period maps 
in the same lap signature class. Then the corresponding global attractors are connection equivalent,
\begin{equation} \label{eq105}
\cA_f\sim\cA_g \ . 
\end{equation} 
\end{theorem}

\bigskip

Here \eqref{eq105} means connection equivalence between global attractors, \cite{firo96,firo00,roc12}. 
This result is not immediate because all the critical elements in $\cP$ fail the hyperbolicity 
condition (H) which would entail the automatic transversality of stable and unstable manifolds 
and the strong Morse-Smale property. Instead, we will address the automatic transversality of 
center stable and center unstable manifolds in the next Sect. 2. 
Then, Theorem \ref{th1} will follow from:

\bigskip

\begin{theorem} \label{th2}
If $f=f(u,u_x)\in\cR$ and $g=g(u,u_x)\in\cR$ belong to the same lap signature class (see 
\cite{roc12,firowo12a}), then there exists a global collective homotopy in $\cR$ between 
$f$ and $g$ which preserves hyperbolicity (H) of all the homogeneous equilibria and 
(normally hyperbolic) frozen rotating waves.   
\end{theorem}

\bigskip

The construction of this {\em global collective homotopy} is very delicate and needs some 
clarification in what concerns the proof in \cite{roc12}, which is restricted to nonlinearities 
of simple type. In Sect. 3 we change and refine some aspects of the proof of Theorem 3 in 
\cite{roc12} making it simpler and amenable to a generalization to nonlinearities of non-simple type. 

In the final Sect. 4 we discuss these results and conclude with the presentation of bifurcation 
diagrams and connection graphs of global attractors for dynamical systems generated by \eqref{eq101}. 
We also describe bifurcation diagrams which reduce a global attractor $\cA_1$ to the trivial global 
attractor $\cA_0=\{0\}$ by a bifurcation homotopy $f^s(u,u_x)=f(s,u,u_x)\in\cR$ with bifurcation 
parameter $1\ge s\ge 0$.

\smallskip

\section{Transversality between center stable and center unstable manifolds}

As already mentioned, the result of Theorem \ref{th1} is not immediate since all the critical elements 
in $\cP$ fail the hyperbolicity condition (H). 
Hence, we replace stable and unstable manifolds by center stable and center unstable manifolds, and 
prove their automatic transversality following closely the previous transversality results established 
for the general problem \eqref{eq101}, see \cite{firowo04,czro08}. See also \cite{hen85,ang86,chchha92} 
for the original transversality results in the case of separated boundary conditions.

We recall that every $2\pi$-periodic solution in $\cP$ is a nonhomogeneous stationary wave frozen 
in time. So, let $\cF$ denote the set of nonhomogeneous stationary waves, and $\cH$ the set of 
heteroclinic orbits between equilibria, either homogeneous or nonhomogeneous. 
Then, for the restrictive set of reversible nonlinearities $f\in\cR$, the global attractor 
$\cA_f$ of \eqref{eq101} decomposes as
\begin{equation} \label{eq201}
\cA_f = \cE\cup\cF\cup\cH \ .
\end{equation}

By reversibility, each nonhomogeneous stationary solution $v(\cdot)\in\cF$ generates a continuum of 
shifted copies $v(\cdot+\vartheta), \vartheta\in S^1$, all normally hyperbolic. This implies that each 
$v(\cdot)\in\cF$ has a unique center manifold composed by its frozen shifted copies. In view of the 
gradient flow variational character of \eqref{eq101} for $f\in\cR$, (see \cite{firaro14}), 
non-stationary orbits $u(t,\cdot)$ converge either to equilibria or stationary waves as 
$t\rightarrow\pm\infty$.
Hence, these orbit connections are essentially heteroclinic connections between equilibria. 
Moreover, the proof of this transversality was already sketched in Propositions 3.1 and 3.2 of 
\cite{firowo04}. In the following, for each $v(\cdot)\in\cF$, we let 
\begin{equation} \label{eq202}
W^{cu}(v(\cdot)) = \bigcup_{\vartheta\in S^1} W^u(v(\cdot+\vartheta)) \quad , \quad 
W^{cs}(v(\cdot)) = \bigcup_{\vartheta\in S^1} W^s(v(\cdot+\vartheta)) \ .
\end{equation}
Finally, we recall that the existence of heteroclinic orbits between two critical elements in 
$\cE\cup\cF$ is determined by a convenient notion of {\em adjacency} between these elements, see 
Theorems 1.3 and 1.4 of \cite{firowo04}.

\bigskip

{\em Proof of Theorem \ref{th1}:}
In the following, we consider only the case of orbit connections between nonhomogeneous stationary 
waves since the remaining cases are simpler. 

We say that two (different) stationary waves $v_\pm(\cdot)\in\cF$ are connected by a heteroclinic 
orbit $u(t,\cdot)$ if this connecting orbit converges to suitably phase shifted stationary waves 
$v^\infty_\pm := v_\pm(\cdot+\vartheta_\pm)$ as $t\rightarrow\pm\infty$ for some fixed 
$\vartheta_\pm\in\R$. 

\smallskip

We introduce here the following definition:

\smallskip

{\bf (D)}: We say that $W^{cu}(v_-)$ and $W^{cs}(v_+)$ are transverse, 
\begin{equation} \label{eq203}
W^{cu}(v_-) \ \otto \ W^{cs}(v_+) \ , 
\end{equation}
if their intersection is empty, or if the strongly unstable manifold $W^u(v^\infty_-)$ and 
the strongly stable manifold $W^s(v^\infty_+)$ are {\em codimension one transverse}, i.e. at 
intersection points their tangential spaces span a codimension one subspace $X_1\subset X$.

\bigskip

Let $u_0=u(0,\cdot)$ denote the initial condition of the connecting orbit. Then, if 
$\cS(t):X\rightarrow X$ denotes the semiflow generated by \eqref{eq101}, we have that 
$u(t,x)=\cS(t)u_0(x)$. 
Moreover, if $D\cS(t)$ is the Fréchet derivative of $\cS(t)$, then for any $w_0\in X$ the curve 
$w(t,\cdot) = \bigl(D\cS(t)u_0(\cdot)\bigr)w_0(\cdot)$ in $X$ defines the classical solution $w(t,x)$ of 
the linearized equation 
\begin{equation} \label{eq204}
w_t = w_{xx} + f_p(u,u_x) w_x + f_u(u,u_x) w \ , \quad x\in S^1 \ , \ t>0 \ , \quad w(0,x) = w_0(x) \ .
\end{equation}
In abstract form this is a linear evolution equation $\dot w = A(t)w$ which generates a solution 
operator $w=\cT(t,\tau)w_\tau, t\ge\tau$. Note that $\cT(t,\tau)$ is injective 
(see for example \cite{czro08}).

We recall that, for the solution $w(t,\cdot)$ of \eqref{eq204}, we have the well known monotone 
nonincreasing property of the zero number 
\begin{equation} \label{eq205}
t \mapsto z(w(t,\cdot)) \ , \ \mbox{is monotone nonincreasing} \ ,
\end{equation} 
(see Lemma 2.6 of \cite{mat82}). Moreover, the zero number strictly decreases at multiple zeros 
of $w(t,\cdot)$, (see \cite{ang88}). 

For $u=v^\infty_\pm$ the equation \eqref{eq204} denotes the linarization around the equilibria 
$v^\infty_\pm\in\cF$ and is autonomous. The corresponding autonomous linear evolution equations with 
linear operators $A^\pm$ generate semigroups $\cT^\pm(t), t\ge 0$, which are analytic. 
We next collect some information regarding the spectral properties of $A^\pm$.

The spectrum of $A^\pm$ is the set of eigenvalues of the second order differential eigenvalue problem 
with periodic boundary conditions
\begin{equation} \label{eq206}
w_{xx} + f_p(u,u_x) w_x + f_u(u,u_x) w = \lambda w \ , \ x\in S^1 \ , 
\end{equation} 
with $u=v^\infty_\pm\in\cF$. 
Let $\spec(A^\pm)=\{\lambda^\pm_j\}_{j=0}^\infty$ denote the set of eigenvalues of \eqref{eq206} 
numbered according to $\lambda^\pm_0>\Re\lambda^\pm_1\ge\Re\lambda^\pm_2\ge\dots$.
By standard spectral theory for \eqref{eq206}, see for example \cite{cole55}, the eigenvalue sequences 
are partially ordered by 
\begin{equation}\label{eq207} 
\Re\lambda^\pm_{2j}>\Re\lambda^\pm_{2j+1} \ , \quad j=0,1,2,\dots \ . 
\end{equation}
As a simple remark we point out that, if $\Re\lambda^\pm_{2j-1}>\Re\lambda^\pm_{2j}$ for any $j\ge 1$, 
then the eigenvalues $\{\lambda^\pm_{2j-1},\lambda^\pm_{2j}\}$ are real.

Let $E^\pm_0$ denote the eigenspaces of constant functions corresponding to $\lambda^\pm_0$, and let 
$E^\pm_j$ denote the generalized eigenspaces corresponding to the spectral sets 
$\{\lambda^\pm_{2j-1},\lambda^\pm_{2j}\}$ for $j\ge 1$. Then, each $w\in E^\pm_j\setminus\{0\}$ has 
only simple zeros on $S^1$ and $z(w)=2j$, \cite{anfi88,mana97}.

\bigskip

Let $u_0\in W^u(v^\infty_-)\cap W^s(v^\infty_+)$. Then $u(t,\cdot)\rightarrow v^\infty_\pm$ as 
$t\rightarrow\pm\infty$ and the linearization \eqref{eq204} around the equilibria $v^\infty_\pm\in\cF$ 
correspond to the asymptotic behaviors of \eqref{eq204} around $u(t,\cdot)$. 

Let $T_{u_0}W^s(v^\infty_\pm)$ and $T_{u_0}W^u(v^\infty_\pm)$ denote the tangent manifolds at $u_0$ of 
the stable and unstable manifolds of $v^\infty_\pm$. By the $S^1$-equivariance of \eqref{eq101} 
these linear manifolds are subspaces of a codimension one subspace $X_1\subset X$.
In fact, for $\vartheta\in S^1$ let $R_\vartheta$ denote the the $S^1$-action induced by the $x$-shift 
$(R_\vartheta u)(x)=u(x+\vartheta)$. Hence, we have
\begin{equation} \label{eq208}
\frac{d}{d\vartheta} (R_\vartheta u_0)(x)|_{\vartheta=0} = 
\frac{d}{d\vartheta} u_0(x+\vartheta)|_{\vartheta=0} = (u_0)_x(x) \ .
\end{equation}
Let $X_0\subset X$ denote the one dimensional linear subspace generated by \eqref{eq208}, 
$X_0 = \span\{(u_0)_x\}$. 
Since $(u_0)_x \not\in (T_{u_0}W^s(v^\infty_\pm) \cup T_{u_0}W^u(v^\infty_\pm))$, this implies
\begin{equation} \label{eq209}
X_0 \cap T_{u_0}W^s(v^\infty_\pm) = X_0 \cap T_{u_0}W^u(v^\infty_\pm) = \{0\} \ .
\end{equation}
Moreover, $T_{u_0}W^s(v^\infty_\pm)$ and $T_{u_0}W^u(v^\infty_\pm)$ span a codimension one 
subspace $X_1$. Since $\dim X_0 = 1$, we obtain
\begin{equation} \label{eq210}
X_0 + X_1 = X \ .
\end{equation}

This also holds for each equilibrium $v\in\cF$ with $X_0=\span\{v_x\}$, and in particular for the 
equilibria $v^\infty_\pm$.

\smallskip

By normal hyperbolicity of the nonhomogeneous stationary wave $v\in\cF$, we define its Morse index 
$i(v)$ as the number $m$ of eigenvalues with strictly positive real part $\Re\lambda_j>0$. 
Then, for $i(v)=m$, $\lambda_{m+1}=0$ corresponds to the eigenfunction $v_{m+1}=v_x(x)$, and we obtain: 
\begin{equation} \label{eq211}
\dim W^{cu}(v) = i(v)+1 \quad , \quad \codim W^{cs}(v) = i(v) \ .
\end{equation}

As in \cite{czro08} (for periodic orbits) two cases need to be considered. 
If the Morse index is odd, $i(v)=2N-1$, then the eigenspace $E_N=E^\pm_N$ corresponds to the 
spectral set $\{\lambda^\pm_{2N-1}=0,\lambda^\pm_{2N}\}$. 
If on the other hand the Morse index is even, $i(v)=2N$, then the eigenspace $E_N$ corresponds to the 
spectral set $\{\lambda^\pm_{2N-1},\lambda^\pm_{2N}=0\}$.
In each case, let $P_0$ denote the projection onto the eigenspace corresponding to the zero eigenvalue, 
i.e. $P_0X=\span\{v_x\}\subset E_N$. Then, for each $w\in P_0X\setminus\{0\}$ we have $z(w)=2N$ and, 
by the monotone nonincreasing property of the zero number \eqref{eq205}, we obtain 
\begin{equation} \label{eq212}
\begin{aligned}
z(v_0) \le &
\begin{cases}
2N-2 = i(v)-1 \quad & \text{if} \ i(v)=2N-1 \ , \\
2N = i(v) \quad & \text{if} \ i(v)=2N \ ,
\end{cases}
\quad  \mbox{for} \ v_0 \in W^u(v) \ , 
\\ 
z(v_0) \ge &
\begin{cases}
2N = i(v)+1 \quad & \text{if} \ i(v)=2N-1 \ , \\
2N+2 = i(v)+2 \quad & \text{if} \ i(v)=2N \ ,
\end{cases}
\quad  \mbox{for} \ v_0 \in W^s(v) \ .
\end{aligned}\end{equation} 

Now assume the existence of an intersection point $u_0$ between the unstable manifold of $v^\infty_-$ 
and the stable manifold of $v^\infty_+$, 
\begin{equation} \label{eq213}
u_0 \in W^u(v^\infty_-) \cap W^s(v^\infty_+) \subset W^{cu}(v_-) \cap W^{cs}(v_+) \ .
\end{equation}
Since the connecting orbit $u(t,\cdot)$ with $u(0,\cdot)=u_0$ lies in 
$W^u(v^\infty_-) \cap W^s(v^\infty_+)$, the invariance of these manifolds implies that 
\begin{equation} \label{eq214}
u_t(t,\cdot) \in T_{u_0}W^u(v^\infty_-) \cap T_{u_0}W^s(v^\infty_+) \ .
\end{equation}
Moreover, since $u_t(0,\cdot)$ is nonzero, the monotone nonincreasing property \eqref{eq205} of the 
zero number implies 
\begin{equation} \label{eq215}
2i(v^\infty_+) < z(u_t(0,\cdot) \le 2i(v^\infty_-) \ .
\end{equation}

Then, combining the inequalities \eqref{eq212} we obtain the following result:

\smallskip

\begin{lemma} \label{le1}
If $u_0\in (W^u(v^\infty_-)\cap (W^s(v^\infty_+)))\setminus \{v^\infty_-,v^\infty_+\}$ and 
$2N^\pm=z(v^\infty_\pm)$, then 
\begin{equation} \label{eq216}
N^-\ge N^+ \quad \mbox{and} \quad i(v^\infty_-) \ge i(v^\infty_+)+1 \ .
\end{equation}
Moreover, if $i(v^\infty_+)=2N^+$, then $N^-\ge N^++1$.
\end{lemma}

Notice that \eqref{eq216} prevents the existence of homoclinic orbits to the set of 
frozen shifted copies of equilibria in $\cF$. 

\bigskip

For the proof of transversality we adapt the proof in \cite{czro08} and use \cite{chchha92} 
(in particular its Theorem 3.1 and the results of Appendix B). 

\smallskip

For any $\psi\in X$ we define 
\begin{equation} \label{eq217}
\nu_\infty(\psi) := \limsup_{n\rightarrow \infty} \Re \log \norm{\cT(n,0)\psi}^\frac1{n} \ .
\end{equation}
The number $\nu_\infty(\psi)$ denotes the {\em Lyapunov characteristic number} for the equation 
\eqref{eq204} with $u=u(t,\cdot)$. 

For any integer $j\ge 0$ we define the spaces 
\begin{equation} \label{eq218}
F^+_j := \{\psi\in X: \nu_\infty(\psi)\le r_j\} \ ,
\end{equation}
where, to respect the spectral set pairs $\{\lambda^+_{2j-1},\lambda^+_{2j}\}$ in $\sigma(A^+)$, 
we let $r_0=\lambda^+_0$ and $r_j=\Re{\lambda^+_{2j-1}}, j\ge 1$. Then, we have 
\begin{equation} \label{eq219}
X = F^+_0 \supset F^+_1 \supset F^+_2 \supset \dots \quad , 
\quad \bigcap_{j=0}^\infty F_j^+ = \{0\} \ .
\end{equation}
Assuming that $\nu_\infty(\psi)\in\{r_j,r_{j+1}\}$ with odd $j$, then 
$\psi\in F^+_j\setminus F^+_{j+1}$ and we consider the asymptotic behavior of 
$v(n,\psi)=\cT(n,0)\psi, n\in\N$, obtaining for some subsequence $k_n\in\N$
\begin{equation} \label{eq220}
\lim_{n\rightarrow\infty} \frac{v(k_n,\psi)}{\norm{v(k_n,\psi)}} = \phi \ ,
\end{equation}
with $\phi$ in the norm one sphere in $X$. This shows that 
\begin{equation} \label{eq221}
z(\psi) \ge z(v(k_n,\psi)) \ge z\Big(\frac{v(k_n,\psi)}{\norm{v(k_n,\psi)}}\Big) = z(\phi) \ .
\end{equation}

Next we use the asymptotic behavior of $v(k_n,\psi)$ to characterize the zero number $z(\psi)$ 
of the initial condition $\psi$.
We first address the zero number characterization of solutions on the stable manifold of the target 
equilibrium $v^\infty_+$.

Again we need to consider the alternatives for the even/odd parity of the index $j$ of the 
eigenvalues $\lambda_j, j\ge 1$.
If $j$ is even, i.e. the eigenspace $E^+_{\frac{j}{2}}$ corresponds to the 
spectral set $\{\lambda^+_{j-1},\lambda^+_j$\}, then 
$\psi\in F^+_{\frac{j}{2}}\setminus F^+_{\frac{j}{2}+1}$ and there exists an asymptotic 
$\phi\in E^+_{\frac{j}{2}}$ such that $z(\psi)\ge z(\phi)=j$. 
If, on the other hand, $j$ is odd, i.e. the eigenspace $E^+_{\frac{j+1}{2}}$ 
corresponds to the spectral set $\{\lambda^+_j,\lambda^+_{j+1}\}$, then 
$\psi\in F^+_{\frac{j+1}{2}}\setminus F^+_{\frac{j+1}{2}+1}$ and there exists an 
asymptotic $\phi\in E^+_{\frac{j+1}{2}}$ such that $z(\psi)\ge z(\phi)=j+1$.

Note that $T_{u_0}W^s(v^\infty_+) = \{\psi\in X: \nu_\infty(\psi)<0\}$.
To continue, let $i(v^\infty_+)=2N^+$ be even. Then $r_{N^+}=\lambda^+_{2N^+}=0$ and for 
a sufficiently large $m\in\N$ we have $T_{u(m,u_0)}W^s(v^\infty_+) = F^+_{N^+}$, 
which is isomorphic to 
\begin{equation}
\cl_X \big(E^+_{N^+} \oplus E^+_{N^++1} \oplus \dots\big) \ .
\end{equation}
In this case we conclude that $z(\psi)\ge 2N^++2$ for $v\in T_{u(m,u_0)}W^s(v^\infty_+)\setminus\{0\}$ 
and $\codim T_{u(m,u_0)}W^s(v^\infty_+) = 2N^++1$.

Before considering the alternative of odd $i(v^\infty_+)=2N^+-1$ we define 
\begin{equation} \label{eq223}
F^+_{\lambda_{2N^+}} := \{\psi\in X: \nu^\infty (\psi)<0\} \ .
\end{equation}
Hence, if $i(v^\infty_+)=2N^+-1$ is odd, then $r_{N^+}<\lambda^+_{2N^+-1}=0$ and for a 
sufficiently large $m\in\N$ we have $T_{u(m,u_0)}W^s(v^\infty_+) = F^+_{\lambda_{2N^+}}$, 
which is isomorphic to 
\begin{equation}
\cl_X \big(\span\{w^+_{2N^++1}\} \oplus E^+_{N^+} \oplus E^+_{N^++1} 
\oplus \dots\big) \ . 
\end{equation}
We conclude that $z(\psi)\ge 2N^+$ for $v\in T_{u(m,u_0)}W^s(v^\infty_+)\setminus\{0\}$ 
and $\codim T_{u(m,u_0)}W^s(v^\infty_+) = 2N^+$.
Since $T_{u_0}W^s(v^\infty_+)$ is the preimage of $T_{u(m,u_0)}W^s(v^\infty_+)$ under the 
injective semiflow $\cT(m,0)$ we obtain the following result:
\begin{lemma}\label{le2}
Assume $\psi\in T_{u_0}W^s(v^\infty_+)\setminus\{0\}$. Then we have the following alternative:  
\begin{equation}\label{eq225}
\begin{cases}
\begin{aligned}
z(\psi) \ge & 2N^++2 \ \mbox{and} \ \codim T_{u_0}W^s(v^\infty_+) =2N^++1 \quad 
\mbox{if} \quad i(v^\infty_+) = 2N^+ \ ,  \\ 
z(\psi) \ge & 2N^+ \ \mbox{and} \ \codim T_{u_0}W^s(v^\infty_+) =2N^+ \quad 
\mbox{if} \quad i(v^\infty_+) = 2N^+-1 \ .
\end{aligned}
\end{cases}
\end{equation}
\end{lemma}

Finally we address the zero number characterization of solutions on the unstable manifold 
of the source equilibrium $v^\infty_-$. 
We invoke the backward unique continuation of the semiflow $\cT(m,0), m\in\N$, and use the 
semiflow $\cT(-m,0)$, which is well defined on the unstable manifold $W^u(v^\infty_-)$. 
Let $i(v^\infty_+)=2N^+$ be even. Then, by Lemma \ref{le1} we have $N^-\ge N^++1$ and, 
for sufficiently large $m\in\N$, there is a subspace $W_{-m}\subset T_{u(-m,u_0)}W^u(v^\infty_-)$ 
such that $z(v)\le 2N^+$ for $v\in W_{-m}$. This implies that, for $W_m:=\cT(m,-m) W_{-m}$ we have 
$z(v)\le 2N^+$ for  $v\in W_m\setminus\{0\}$. In view of Lemma \ref{le2}, this shows that 
$W_m\cap W^s(v^\infty_+)=\{0\}$ which implies $W_m\subseteq T_{u(m,u_0)}W^u(v^\infty_+)$. 
Therefore, for $m\in\N$ sufficiently large, $W_m$ is isomorphic to 
\begin{equation}
E^+_0 \oplus \dots \oplus E^+_{N^+-1} \oplus \span\{w^+_{2N^+-1}\} \ ,
\end{equation}
where $w^+_{2N^+-1}$ is an eigenfunction associated to the (real) eigenfunction $\lambda^+_{2N^+-1}$. 
This shows that
\begin{equation}
\dim W_{-m} = \dim W_m = 2N^+ \ , \quad \mbox{and} \quad z(v)\le 2N^+ \ , \ 
v\in W_{-m}\setminus\{0\} \ .
\end{equation}
Alternatively, let $i(v^\infty_+)=2N^+-1$ be odd. By \eqref{eq216} we have $N^-\ge N^+$ and, 
as in the previous case, for sufficiently large $m\in\N$ there is a subspace 
$W_{-m}\subset T_{u(-m,u_0)}W^u(v^\infty_-)$ such that $z(v)\le 2N^+-2$ for $v\in W_{-m}$. 
This shows that $z(v)\le 2N^+-2$ for  $v\in W_m\setminus\{0\}$, where $W_m:=\cT(m,-m) W_{-m}$ is 
isomorphic to
\begin{equation}
E^+_0 \oplus \dots \oplus E^+_{N^+-1} \ .
\end{equation}
This implies that
\begin{equation}
\dim W_{-m} = \dim W_m = 2N^+-1 \ , \quad \mbox{and} \quad z(v)\le 2N^+-2 \ , \ 
v\in W_{-m}\setminus\{0\} \ .
\end{equation}

In summary, we obtained the following result: 
\begin{lemma}\label{le3}
There is a subspace $W_0=\cT(0,-m)W_{-m}\subset T_{u_0}W^u(v^\infty_-)$ with $\dim W_0 = i(v^\infty_+)$ 
such that, if $\psi\in W_0\setminus\{0\}$, then the following alternative holds: 
\begin{equation}\label{eq230}
\begin{cases}\begin{aligned}
z(\psi) \le & 2N^+ \quad \mbox{if} \quad i(v^\infty_+) = 2N^+ \ ,  \\ 
z(\psi) \le & 2N^+-2 \quad \mbox{if} \quad i(v^\infty_+) = 2N^+-1 \ . 
\end{aligned}
\end{cases}
\end{equation}
\end{lemma}

\smallskip  

The automatic transversality of center stable and center unstable manifolds, claimed in the next 
Theorem \ref{th3}, is now obtained from the combination of Lemmas \ref{le2} and \ref{le3}.
Note that this transversality implies the preservation of heteroclinic orbit connections in the 
$S^1$-equivariant case and entails the connection  equivalence \eqref{eq105} asserted by 
Theorem \ref{th1}. 

\smallskip

\begin{theorem} \label{th3}
Let $v_\pm\in\cF$ denote normally hyperbolic fixed points of the semigroup $\cS(t)$ generated by 
\eqref{eq101}. Then the center unstable manifold of $v_-$, $W^{cu}(v_-)$, and the center stable 
manifold of $v_+$, $W^{cs}(v_+)$, intersect transversely:
\begin{equation} \label{eq231}
W^{cs}(v_+) \ \otto \ W^{cu}(v_-) \ .
\end{equation} 
\end{theorem}

\smallskip

{\em Proof of Theorem \ref{th3}:}
Recall that if $W^{cu}(v-) \cap W^{cs}(v_+) = \emptyset$ then $W^{cu}(v_-)$ and $W^{cs}(v_+)$ are 
transverse by definition (D). 
So, we assume the existence of an intersection point $u_0$ between the center unstable manifold of $v_-$ 
and the center stable manifold of $v_+$.

Then, by \eqref{eq230} for $w\in W_0\setminus\{0\}\subset T_{u_0}W^u(v^\infty_-)\setminus\{0\}$ and 
\eqref{eq225} for $w\in T_{u_0}W^s(v^\infty_+)\setminus\{0\}$ we obtain
\begin{equation} \label{eq232}
W_0 \cap T_{u_0}W^s(v^\infty_+) = \{0\} \ .
\end{equation}
Moreover, the combination of \eqref{eq225} and \eqref{eq230} imply
\begin{equation} \label{eq233}
\dim W_0 = \codim T_{u_0}W^s(v^\infty_+)-1 \ .
\end{equation}
This shows that 
\begin{equation} \label{eq234}
W_0 \oplus T_{u_0}W^s(v^\infty_+) = \widetilde X_1 \ ,
\end{equation}
where $\widetilde X_1\subset X$ has $\codim \widetilde X_1=1$. Since $W_0\subset 
T_{u_0}W^u(v^\infty_-)$, this implies $W^u(v^\infty_-) \ \otto \ W^s(v^\infty_+)$ and, by our 
definition (D), this shows the transversality result \eqref{eq231} and completes the proof of 
Theorem \ref{th3}. 

To finish, let $\widetilde X_0$ denote the one dimensional complement of $\widetilde X_1$, i.e. 
$\widetilde X_0 + \widetilde X_1 = X$, obtained by using again the backward linear flow generated by 
\eqref{eq204} with $u=u(t,\cdot)$. Then, we have
\begin{equation} \label{eq235}
T_{u_0}W^u(v^\infty_-) + T_{u_0}W^s(v^\infty_+) + \widetilde X_0 = X \ ,
\end{equation} 
concluding this discussion of transversality.
$\hfill \square$

\bigskip

Notice that forty years ago Dan Henry \cite{hen85} used the expression ``amazing'' to refer to this 
type of result, and advised the reader to ``look at it again''!

\smallskip

\section{Homotopy construction}

We first consider the obstructions to the homotopy construction between a reversible nonlinearity 
$f_0=f_0(u,u_x)$ and an $f_1=f_1(u)$ in the class of Hamiltonian systems. 

Let $T=T(u_0):\cD\rightarrow\R_+$ denote the period map of a reversible nonlinearity $f=f(u,u_x)$.
By continuity, we extend $\cD$ to $u_0=e$ corresponding to the center $(e,0)$ which is encircled by 
periodic orbits. For simplicity, by a translation $u\mapsto u-e$ we let $e=0$. 

A necessary and sufficient condition for $T=T(u_0)$ to be realizable by a Hamiltonian nonlinearity 
$g=g(u)$ is

\bigskip

{\bf (C)}: \quad $d(u_0T(u_0))/du_0>0, u_0\in\cD$. 

\bigskip

For references see Theorem 4.2.5 and Proposition 4.2.6 of \cite{sch90}, Theorem 5 of \cite{ura64}, 
Theorem 4.2 of \cite{furo91}, and Section 4 of \cite{roc07}.

\bigskip

We first show the existence of reversible nonlinearities $f=f(u,u_x)$ for which this condition is not 
satisfied. For this purpose, consider the nonlinearity 
\begin{equation} \label{eq301}
f(u,u_x)=\frac{k}{2}u\left(ku^2+\sqrt{k^2u^4+4u_x^2}\right) \ . 
\end{equation} 
Each level curve of the first integral $I(u,v)=\frac1{2}\left(ku^2+\sqrt{k^2u^4+4v^2}\right)$ of 
\eqref{eq102} is an ellipse $v^2+kI_0u^2=I_0^2$ where $I_0:=I(u_0,0)=ku_0^2$. Therefore, the 
period map for this nonlinearity is given by 
\begin{equation} \label{eq302}
T(u_0) = 4 \int_0^{u_0} \frac{du}{k\sqrt{u_0^4-u_0^2u^2}} = \frac{2\pi}{k|u_0|} \ , \quad u_0\ne 0 \ .
\end{equation} 
Hence, since $d(u_0T(u_0))/du_0=0$ for $u_0\ne 0$, condition (C) is not satisfied. 

\bigskip

If the period map $T(u_0)$ satisfies condition (C) we are able to find a Hamiltonian pendulum 
nonlinearity $g=g(u)$ with the same period map, \cite{roc07}. See also \cite{ura64,sch90} for the 
preliminary results. In this case, the convex combination 
\begin{equation} \label{eq303}
f_\tau(u,u_x)= (1-\tau)f(u,u_x)+\tau g(u)
\end{equation} 
provides the desired {\em pendulum realization homotopy} in the Hamiltonian class of nonlinearities. 

On the other hand, the failing of condition (C) prevents the homotopy in the Hamiltonian 
class since the period map is not realizable in this class. Therefore, the homotopy has to be 
constructed directly in $\cR$ using the phase portrait of \eqref{eq102}. We define the 
{\em Hamiltonian realization homotopy} as a continuous transformation of the cyclicity set $\cC$ of 
$f(u,u_x)$, which preserves the hyperbolicity of all homogeneous equilibria and frozen rotating waves, 
and shrinks the period map along the $u$-axis until it satisfies condition (C).

\bigskip

Let $K=(n-1)/2+r$ denote the number of connected regions of the cyclicity set $\cC$. Here $n$ denotes 
the odd number of equilibria and $r$ the number of annular regions. Annular regions are the multiply 
connected regions bounded by two saddles and their homoclinic ODE separatrices (see for example the 
region $\cC_1$ in Fig.6).
Notice that each equilibrium of \eqref{eq101} corresponds to a center or a saddle point of the phase 
portrait of \eqref{eq102}. The centers and saddles alternate for increasing values of $u_0$, starting 
and ending with saddles. 
Then $(n-1)/2$ is the number of punctured disks encircling the centers. Moreover, $r$ is the number 
of annular regions which surround more than one single center. This implies that 
$0 \le r \le (n-3)/2$. See \cite{firowo12a}. 

Finally, a reversible nonlinearity $f=f(u,u_x)$ is of {\em simple type} if $r=0$, i.e. if there are 
no annular regions. This is a slightly more restrictive definition then the one used in \cite{roc12}. 
If $r=0$ then each $2\pi$-periodic orbit of \eqref{eq102} encircles exactly one center in the phase 
plane $(v,v_x)$. 

Then, let $g=g(u)$ denote a Hamiltonian nonlinearity with period map $T=T(u_0)$, and let $G$ denote 
the potential function of \eqref{eq102}, i.e. $G'=g$. Under the hyperbolicity assumption (H) for 
\eqref{eq101} with $f=f_1(u,u_x)$, all the zeros of $g$ are nondegenerate (i.e. $g(u_0)=0$ implies 
$g'(u_0)\ne 0$).
Moreover, the period map $T=T(u_0)$ satisfies the nondegeneracy condition: $T(u_0)=\frac{2\pi}{k}$ 
implies $T'(u_0)\ne 0$ for all positive integers $k\in\N$. See \cite{roc07} for details. 

For simplicity, we make here the generic assumption 

\bigskip

{\bf (M)}: the potential function $G$ is a {\em Morse function}. 

\bigskip

Hence, in addition to the nondegeneracy of the zeros of $g$, we also assume that all the critical values 
of $G$ are distinct ($g(u_0)=g(u_1)=0$ implies $G(u_0)\ne G(u_1)$). Since we can always add a small 
perturbation to $f_1$, this generic assumption does not affect our result.

Initially, we consider that the lap signature class of $f(u,u_x)$ is of simple type. This implies that 
$r=0$ and we deal with a cyclicity set composed only of isolated punctured disks around the centers 
in the phase portrait.

\bigskip 

\begin{figure}[ht] \label{fig1}
\begin{center}
\includegraphics[scale = .7]{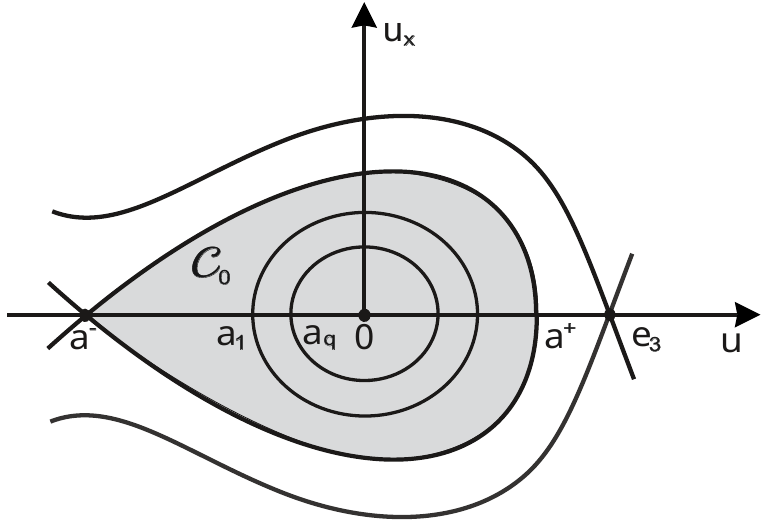}
\end{center}
\caption{\small Cyclicity set $\cC_0$ in the case of $n=3, r=0$ and with $q=2$ $2\pi$-periodic orbits.}
\end{figure}

{\em Proof of Theorem \ref{th2}:} The proof essentially consists on the construction of a Hamiltonian  
realization homotopy in $\cR$ between $f_0(u,u_x)$ and a nonlinearity $f_1(u,u_x)$ satisfying 
condition (C) of a Hamiltonian nonlinearity $g=g(u)$. 
This homotopy preserves: (a) hyperbolicity of all critical elements of the flow generated by 
\eqref{eq101} with $f=f_0(u,u_x)$; and hence, (b) the lap signature of the period map $T=T_f$. 

We consider first the case of a single connected cyclicity set $\cC_0$ surrounding a unique center, 
{\em Case (I)}, and then we consider the case of multiple isolated punctured disks, {\em Case (II)}. 
Later on we will treat the case of nonlinearities with lap signature class of non-simple type.

\bigskip

{\em Case (I): $n=3$ and $r=0$.} For the reversible nonlinearity $f=f_0(u,u_x)$ let 
$e_1 < e_2=0 < e_3$ denote the three zeros of $f_0(\cdot,0)$ corresponding to two saddle points, 
${\bf e}_1, {\bf e}_3$, and a center at the origin ${\bf e}_2=(0,0)$. The boundary $\partial\cC_0$ 
of the cyclicity set is given by a saddle point and an orbit homoclinic to this saddle. 
Without loss of generality we assume that the saddle point is ${\bf e}_1$.
Therefore, the cyclicity set $\cC_0$ is positioned to the right of ${\bf e}_1$, see Fig. 1. 

For simplicity we define the interval $(a^-,a^+)$ corresponding to the $u$-values of the cyclicity 
set $\cC_0$. 
Here $a^-:=e_1$ and $a^+, 0<a^+<e_3$, corresponds to the maximum $u$-value of the homoclinic orbit. 
Then, the period map for $f_0(u,u_x)$ satisfies $T:(a^-,0)\cup(0,a^+)\rightarrow\R$ and we extend 
the domain of $T(u_0)$ to $u_0=0$ by continuity. By the hyperbolicity assumption (H) and the 
smoothness of $f_0(u,u_x)$ we obtain, 
\begin{equation} \label{eq304}
T(0)\in \left(\frac{2\pi}{k},\frac{2\pi}{k-1}\right) \ \mbox{ for some } k\in\N \ , \ 
\mbox{ and } \ T'(0)=0 \ ,
\end{equation} 
see for example \cite{sch90}. 

\begin{figure}[ht] \label{fig2}
\begin{center}
\includegraphics[scale = .9]{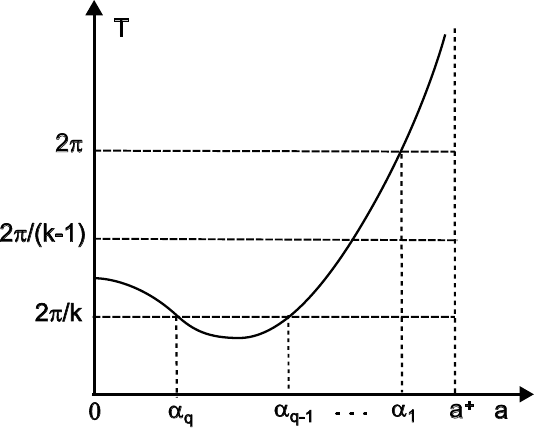}
\end{center}
\caption{\small Graph of the period map $T=T(a)$ on the half interval $(0,a^+)$. 
Here $\alpha_j=\max_{x\in[0,2\pi]}v_j(x,a_j), j=1,\dots,q$, denote the maximal values of the 
$2\pi$-periodic solutions.}
\end{figure}

Assume the set of periodic solutions $\cP$ to be nonempty, i.e. $q\ge 1$.
We recall that the frozen rotating waves of $f_0(u,u_x)$ are the $2\pi$-periodic solutions of 
\eqref{eq102}, ${\bf v}_j\in\cP$ for $1\le j\le q$. Let $a_j=\min_{x\in[0,2\pi]}v_j(x)$ denote the 
minimal initial values of the $u$-coordinates of these $2\pi$-periodic solutions. 
Then, we have $a^-<a_1<\dots<a_q<0$ and for each $1\le j\le q$ we have $T(a_j)=\frac{2\pi}{k_j}$ 
for some $k_j\in\N$. 

The proof of case {\em (I)} then proceeds by the construction of a homotopy between $f_0$ and a 
reversible $f_1$ which preserves hyperbolicity (H) and for which the period map $T=T_1(a)$ satisfies 
condition (C). This involves the use of a local diffeomorphism of the plane $(u,u_x)$ which preserves 
the period map values. 
Based on \cite{firowo04}, Lemma 5.1, this diffeomorphism is described in \cite{roc12}, Sect. 5, and 
is recalled here for the benefit of the reader. 

Let $(v,p(v)) = (v(\cdot,a),p(v(\cdot,a))$ denote the periodic orbits of \eqref{eq102} on the phase 
plane $(u,u_x)$ where $v=v(\cdot,a)$ denotes the solution with $v(0,a)=a, v_x(0,a)>0$. 
On the cyclicity set $\cC_0$ we define the scaling map 
\begin{equation} \label{eq305}
\Phi(v,p) = \Omega(v,p) (v,p) 
\end{equation} 
where $\Omega: \cC_0\rightarrow\R$ is constant along the periodic orbits of \eqref{eq102}, i.e. 
\begin{equation} \label{eq306}
\Omega(v(\cdot,a),p(v(\cdot,a))) = \Omega(a,0) \ . 
\end{equation} 
Then, let $\omega: (0,a^+)\rightarrow\R$ denote the scale function $\omega(a)=\Omega(a,0)$ which is 
assumed monotone nondecreasing in order to have $\Phi$ as a diffeomorphism on the cyclicity set $\cC_0$. 
We extend the domain of $\omega$ to the interval $(a^-,a^+)$ using the minimal values of the periodic 
orbits in the cyclicity set and defining $\omega(0):=\omega(0^+)=\omega(0^-)$. 
Finally, we extend the diffeomorphism $\Phi$ to the whole phase space $(u,u_x)\in\R^2$ by defining 
$\Phi$ as the identity in $\R^2\setminus\cC_0$.

To define the scale function $\omega|_{(0,a^+)}$, let $c_1, c_2\in(\alpha_1,a^+)$ denote two constants 
satisfying $\alpha_1<c_1<c_2<a^+$, where again $\alpha_1$ denotes the maximal value of the outermost 
$2\pi$-periodic solution $v_1(\cdot,a_1)$. 
Recall that $\lim_{a\rightarrow a^+}T(a) = +\infty$, hence $T(a)>2\pi$ for $a\in(\alpha_1,a^+)$, 
see Fig. 2. 
Changing the phase portrait of \eqref{eq102} in a small neighborhood of the homoclinic orbit to 
${\bf e}_1$ does not affect the global attractor $\cA_f$. 
Therefore, we can assume that the period map satisfies $\lim_{a\rightarrow a^+}T(a) = +\infty$ 
monotonically. 
Hence, we have $T'(a)>0$ in a neighborhood of $a^+$ and we always choose $c_1$ in this neighborhood. 
Then, for a parameter $\delta\in(0,1)$, we define 
\begin{equation} \label{eq307}
\omega(a) = \left\{ \begin{array}{ll}
\delta \quad \ & \text{for} \ a\in(0,c_1) \ , \\ \\
\text{$C^2$-smooth monotone increasing} \quad \ & \text{for} \ a\in[c_1,c_2] \ , \\ \\
1 \quad \ & \text{for} \ a\in(c_2,a^+) \ . 
\end{array}\right. 
\end{equation} 

\begin{figure}[ht] \label{fig3}
\begin{center}
\includegraphics[scale = .9]{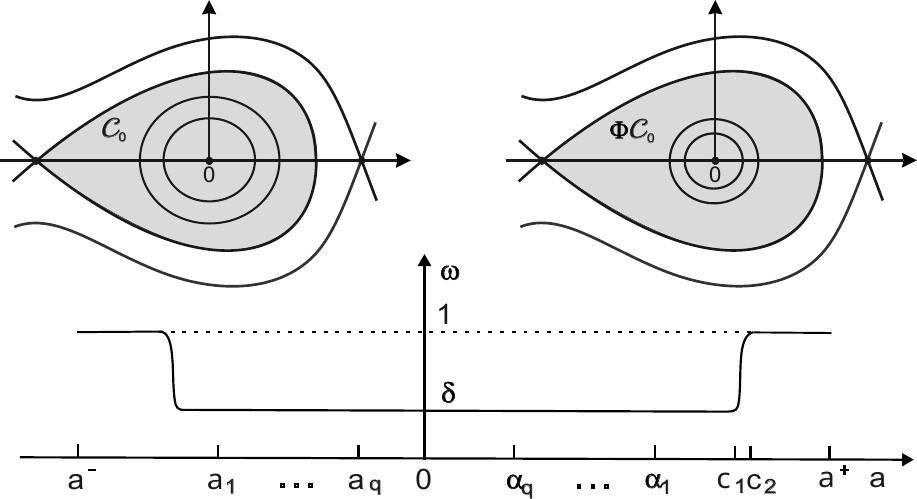}
\end{center}
\caption{\small Upper left: Cyclicity set $\cC_0$ for $f=f_0(u,u_x)$. Upper right: 
Cyclicity set $\Phi\cC_0$ for $f=f_1(u,u_x)$. Bottom: Graph of the scale function $\omega$.}
\end{figure}

Therefore, with this scale function $\omega$, the diffeomorphism $\Phi$ shrinks all the $2\pi$-periodic 
orbits to a neighborhood of the origin as $\delta$ decreases, see Fig. 3. 
Since the scaling of the period map preserves the nondegeneracy of the points corresponding to the 
maximal values of the $2\pi$-periodic solutions, their hyperbolicity is preserved by the scaling. 
Moreover, the period map $T_1(a)$ satisfies 
\begin{equation} \label{eq308}
T_1(a) = T(\omega(a)a) \ .
\end{equation} 
Hence, we obtain
\begin{equation} \label{eq309}
\frac{d}{da} \left(aT_1(a)\right) = T(\omega(a)a) + (\omega(a)+a\omega'(a))aT'(\omega(a)a) \ .
\end{equation}
Notice that $\omega'(a)=0$ for all $a>0$ except $a\in(c_1,c_2)$ where $\omega'(a)>0$ and, by our 
choice of $c_1$ in a neighborhood of $a^+$, the period map satisfies $T'(\omega(a)a)>0$ there.
This implies that, for $\delta$ sufficiently small, $T_1(a)$ satisfies condition (C). 
We conclude that $f_1(u,u_x)$ is realizable in the Hamiltonian class of nonlinearities completing 
the proof of case {\em (I)}.

\bigskip

{\em Case (II): $n\ge 5$ and $r=0$.} Now we deal with a cyclicity set 
\begin{equation} \label{eq310}
\cC=\bigcup_{1\le k\le (n-1)/2} \cC_k
\end{equation}  
composed of several punctured disks surrounding the $(n-1)/2$ centers of \eqref{eq102}, see Fig. 4. 
Since each $\cC_k$ is isolated (by condition (M) and the simple type $r=0$), the proof in this case 
follows by repeating the previous procedure in each region $\cC_k, k=1,\dots,(n-1)/2$. Therefore, we 
obtain $(n-1)/2$ disjoint graphs of Hamiltonian realizations, one for each region $\cC_k$. Then, a 
Hamiltonian $g(u)$ which realizes $f_1(u,u_x)$ is obtained by extending the domain to $\R$ joining 
these $g$ graphs. For this extension, we fill the $(n-3)/2$ closed interval gaps with $C^2$-smooth 
functions without introducing further zeros of $g$. Similarly, in the two remaining unbounded 
intervals, our extension choice is such that no further zeros are introduced. 
Finally, we choose a globally bounded extension $g:\R\rightarrow\R$. 

\begin{figure}[ht] \label{fig4}
\begin{center}
\includegraphics[scale = .9]{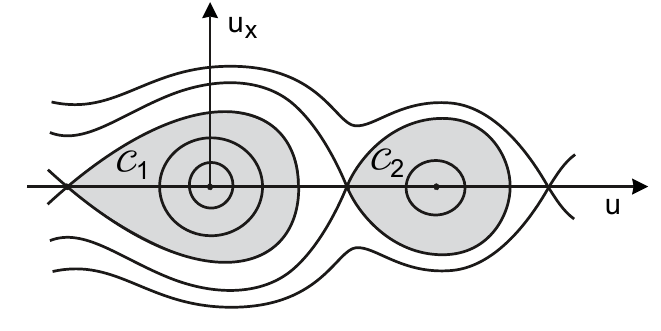}
\end{center}
\caption{\small Cyclicity set $\cC=\cC_1\cup\cC_2$ for $f=f_0(u,u_x)$ with $n=5$ and $r=0$.}
\end{figure}

Clearly, $f_0(u,u_x)$ and $g(u)$ belong to the same lap signature class. 
Then, the desired Hamiltonean realization homotopy has the form (see \cite{roc12}),
\begin{equation} \label{eq311}
f_\tau(v,p) = \Omega_\tau(v,p)\left(f_0\circ\Phi_\tau^{-1}(v,p)\right) \ , \ 0\le\tau\le 1 \ , 
\end{equation} 
where
\begin{equation} \label{eq312}
\Phi_\tau(v,p) = \Omega_\tau(v,p)\,(v,p) \ , \quad \Omega_\tau(a,0) = \omega_\tau(a) \ , 
\quad \omega_\tau = (1-\tau)+\tau\omega \ ,
\end{equation} 
followed by the pendulum realization homotopy \eqref{eq303} between $f_1(u,u_x)$ and $g(u)$.
This completes the proof of case {\em (II)}. 

\bigskip

{\em Remark:} Consider again the Hamiltonian realization of $f=f_0(u,u_x)$ by $g(u)$ in case 
{\em (I)}. 
The cyclicity set $\cC_0$ is {\em right oriented} if its boundary $\partial\cC_0$ contains the saddle 
point ${\bf e}_1$, i.e. $a^-=e_1$. Similarly, $\cC_0$ is {\em left oriented} if $\partial\cC_0$ contains 
the saddle point ${\bf e}_3$, i.e. $a^+=e_3$. 
Then, if $\cC_0$ is right oriented, the nonlinearity $g$ satisfies $g(a^-)=g(e_1)=0$ and is 
negative for $u\in(a^-,0)$ and positive for $u\in(0,a^+)$. Therefore, the extension of $g$ to the 
interval gap $(a^+,e_3)$ is also positive. 
On the other hand, if $\cC_0$ is left oriented, we have $g(e_2)=g(a^+)=0$ and the extension of 
$g$ to the interval gap $(e_1,a^-)$ is negative. 
This holds in the multiple case {\em (II)} for all the isolated components $\cC_k$ of the cyclicity set.

\smallskip

Due to the integrability and the simple type of $f_0$ the phase portrait of \eqref{eq102} for 
$f=f_0(u,u_x)$ has the following characteristic. There is a saddle point ${\bf e}_m, 1\le m\le n$ 
and odd, such that all $\cC_k$ to the left of ${\bf e}_m$ are right oriented, and all $\cC_k$ to the 
right of ${\bf e}_m$ are left oriented. See Fig. 5.

\begin{figure}[ht] \label{fig5}
\begin{center}
\includegraphics[scale = .9]{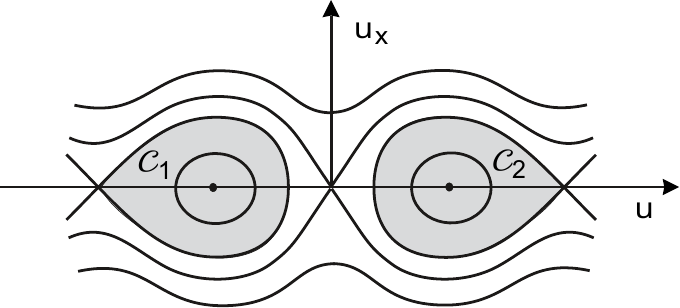}
\end{center}
\caption{\small Cyclicity set $\cC=\cC_1\cup\cC_2$ for $f=f_0(u,u_x)$ with $n=5$, $r=0$ and $m=3$. 
The previous Fig. 4 illustrates the case $m=n=5$.}
\end{figure}

The phase portrait of \eqref{eq102} for $f=g(u)$ depends essentially on the singular values 
of the potential function 
\begin{equation} \label{eq313}
G(a)=\int_{e_2}^a g(s) ds 
\end{equation} 
at the saddle points, i.e. the local maxima of $G$, and not on the Morse type of $G$. Indeed, the 
Morse type of $G$ is easily modified by changing the values of the extended $g$ in its $(n-1)/2$ 
closed interval gaps. 

Since all $\cC_k$ are right oriented at the left of ${\bf e}_m$ and left oriented at the right, all 
the extensions of $g$ are positive to the left of $e_m$ and negative to the right. Then, by 
\eqref{eq313}, ${\bf e}_m$ is the saddle point with the maximum value 
$G(e_m)=\max\{G(e_{2k-1}): 1\le k\le (n+1)/2\}$ and the sequence of singular values $G(e_{2k-1})$ 
satisfies the inequalities 
\begin{equation} \label{eq314}
G(e_1) < \dots < G(e_m) \quad , \quad G(e_m) > \dots > G(e_n) \ .
\end{equation} 

This ensures that the phase portraits of \eqref{eq102} for $f=f_0(u,u_x)$ and $f=g(u)$ are 
qualitatively the same. 

\bigskip

We now prove the generalization of our main result to the case of nonlinearities with lap signature 
of non-simple type. We prove this following the same approach used in the previous cases {\em (I)} 
and {\em (II)}. 
We first consider the case of $\cC$ with a single outermost annular region $\cC_1$ surrounding two 
cyclicity sets $\cC_2\cup\cC_3=\cC\setminus\cC_1$, case {\em (III)}, and then we consider the case of 
multiple outermost annular regions of $\cC$ and punctured disks, case {\em (IV)}.

\smallskip

In view of the proof of previous cases, we only need to consider nonlinearities $f_0(u,u_x)$ of 
non-simple type, $r\ge 1$. This implies $n\ge 5$. As before we construct a Hamiltonian realization 
homotopy in $\cR$ between $f_0(u,u_x)$ and a nonlinearity $f_1(u,u_x)$ satisfying condition (C). 

\bigskip

{\em Case (III): $n=5$ and $r=1$.} In this case, the phase portrait of \eqref{eq102} for 
$f=f_0(u,u_x)$ has three cyclicity regions. In addition to two punctured disks $\cC_2, \cC_3$ 
surrounding the centers ${\bf e}_2, {\bf e}_4$, there is an annular region $\cC_1$ surrounding both 
punctured disks. 
This implies that the saddle point ${\bf e}_3$ has two homoclinic orbits composing the boundary 
$\partial(\cC_2\cup\cC_3)$. The cyclicity regions $\cC_2$ and $\cC_3$ have opposite left/right 
orientations and the homoclinic orbits form an $\infty$ shaped curve. See Fig. 6.
 
\begin{figure}[ht] \label{fig6}
\begin{center}
\includegraphics[scale = 1]{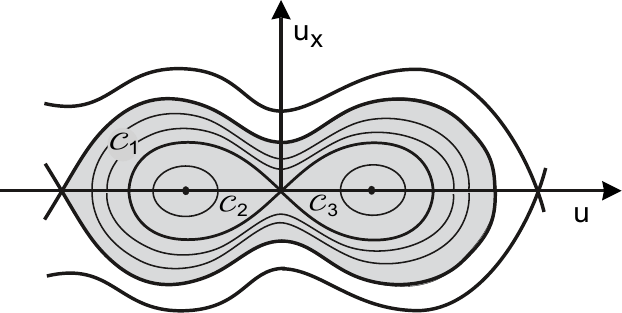}
\end{center}
\caption{\small Cyclicity set $\cC=\cC_1\cup\cC_2\cup\cC_3$ for $f=f_0(u,u_x)$ with $n=5$, $r=1$ and 
$q=4$. Two of these $2\pi$-periodic orbits are in the annular region $\cC_1$.}
\end{figure}

By a translation in $u$, we assume here that $e_3=0$. Therefore, the $n=5$ equilibria of \eqref{eq101} 
for $f=f_0(u,u_x)$ satisfy $e_1<e_2<e_3=0<e_4<e_5$. We use the same notation that was used in the 
previous cases. 
Specifically, we let $(a^-,0)=(e_1,0)$ denote the first saddle point, and $(a^+,0)$ the maximum value 
of the orbit homoclinic to $(e_1,0)$. Moreover, we denote by $a_1, \dots, a_{q_1}$ the Neumann initial 
values of the $2\pi$-periodic orbits $v_j\in\cP$ in the annular region $\cC_1$, i.e. 
$a_j=\min_{x\in[0,2\pi]}v_j, 1\le j\le q_1$. 
Similarly, we denote by $\alpha_1,\dots,\alpha_{q_1}$ the corresponding maximum values 
$\alpha_j=\max_{x\in[0,2\pi]}v_j, 1\le j\le q_1$. Notice that in annular regions $\cC_k$ the numbers 
$q_k$ are always even, see \cite{firowo12a}.

The proof follows the same argument employed in case {\em (I)}. We use again a shrinking scale 
function $\omega$ which grants the Hamiltonian realization of \eqref{eq307} for $f=f_1(u,u_x)$. 
The essential difference here is the shrinking of all equilibria and $2\pi$-periodic orbits in 
$\cC$ to a neighborhood of the middle saddle point ${\bf e}_3$ instead of a neighborhood of the 
center ${\bf e}_2$.

Then, let again $c_1,c_2$, denote the two constants satisfying $\alpha_1<c_1<c_2<a^+$ ($c_1$ in the 
appropriate neighborhood of $a^+$) and, for $\delta\in(0,1)$, define the scale function $\omega$ by 
\eqref{eq307}.
\begin{figure}[ht] \label{fig7}
\begin{center}
\includegraphics[scale = .9]{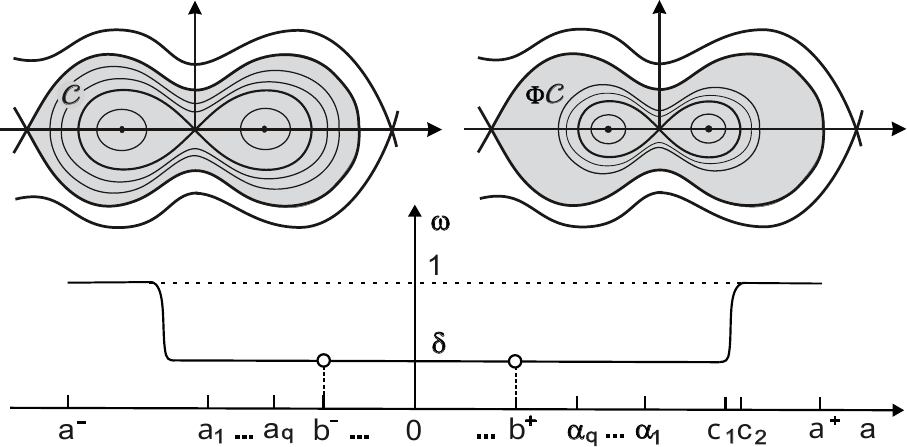}
\end{center}
\caption{\small Upper left: Cyclicity set $\cC$ for $f=f_0(u,u_x)$. Upper right: 
Cyclicity set $\Phi\cC$ for $f=f_1(u,u_x)$. Bottom: Graph of the scale function $\omega$. 
Here $b^-$ and $b^+$ denote respectively the minimum and maximum $u$-values of the orbits homoclinic 
to the origin. The white balls correspond to these homoclinic orbits, and $b^\pm\not\in\cD$ since the 
homoclinic orbits are not in $\cC$. Also note that $0\not\in\cD$.}
\end{figure} 
Once more, using this scale function $\omega$, we define the diffeomorphism $\Phi:\cC\rightarrow\cC$ 
which we extend to $\cl\cC$ by continuity, and to the phase plane $\R^2$ by the identity, 
$\Phi|_{(\R^2\setminus\cl\cC)}=\id$.

As expected, $\Phi$ shrinks all equilibria and $2\pi$-periodic orbits in $\cC$ to a neighborhood 
of the origin $(e_3,0)=(0,0)$ as $\delta$ decreases, preserving hyperbolicity, see Fig. 7. 
In addition, under $\Phi$ the period map $T_1$ obtained from the period map $T$ of \eqref{eq102} 
for $f=f_0(u,u_x)$ satisfies \eqref{eq308}. Then, from \eqref{eq309} we again conclude that $T_1(a)$ 
satisfies condition (C). Hence $f_1(u,u_x)$ is realizable by a Hamiltonian nonlinearity 
$g(u)$ and the desired homotopy from $f_1(u,u_x)$ to $g(u)$ has the form \eqref{eq311}, 
\eqref{eq312}, which is followed by the pendulum realization homotopy \eqref{eq303}. 
This completes the proof in case {\em (III)}.

\bigskip

{\em Case (IV): $n\ge 7$ and $r\ge 1$.} 
We proceed sequentially, following the total order imposed by the regular parenthesis structure of 
the cyclicity regions. 
For example, the regular parenthesis structure of the cyclicity region shown in Fig. 6 is 
\begin{equation} \label{eq315}
(( ) ( )) \ .
\end{equation}  
The use of the regular parenthesis structure of the cyclicity regions instead of the 
$2\pi$-periodic orbits is necessary to overcome the case of $r\ge 1$ without $2\pi$-periodic 
orbits in some annular region $\cC_k$, i.e. $q_k=0$.

We initially consider the first region $\widetilde\cC$ which may contain the first and outermost 
$2\pi$-periodic orbit corresponding to the solution with minimal value 
$a_1=\min_{x\in[0,2\pi]}v_1(x,a_1)$. 
If $\widetilde\cC$ is an isolated punctured disk we apply the scale shrinking procedure described 
in case {\em (II)}. So, next we assume that $\widetilde\cC$ is an annular region.

We let $\cC$ denote the union of $\widetilde\cC$ with all the cyclicity regions that it encloses. 
Then, we apply the same argument used in case {\em (III)} shrinking all equilibria and $2\pi$-periodic 
orbits contained in $\cC$ to a neighborhood of the unique saddle point in the inner boundary of 
$\widetilde\cC$. 
For example, in Fig. 6 the saddle point is ${\bf e}_3$. 
Therefore, using the same notation as in case {\em (III)}, we obtain a diffeomorphism 
$\Phi:\cC\rightarrow\cC$ which preserves hyperbolicity and leads to a period map $T_1|_{(a^-,a^+)}$ 
satisfying condition (C). 

Let $\kappa$ denote the number of annular regions and punctured disks with an outermost homoclinic 
orbit in their boundaries. Let $\widetilde\cC_k$ denote these either punctured disks or outermost 
annular regions of the cyclicity set $\cC$. If $\widetilde\cC_k$ is a punctured disk we define 
$\cC_k:=\widetilde\cC_k$.  
If $\widetilde\cC_k$ is an annular region we define $\cC_k=\widetilde\cC_k\cup\overline\cC_k$ where 
$\overline\cC_k$ denotes the union of all the cyclicity sets encircled by $\widetilde\cC_k$. In this 
way we obtain a sequence $\cC_1,\dots,\cC_\kappa$ were we can apply repeatedly the above procedure. 
Therefore, we obtain  
\begin{equation} \label{eq316}
\cC=\bigcup_{1\le k\le\kappa} \cC_k
\end{equation} 
with diffeomorphism $\Phi:\cC\rightarrow\cC$ which preserves hyperbolicity. Moreover, the period maps 
$T_1(a)$ restricted to the intervals defined by the maximal homoclinic $u$-values satisfy condition (C) 
for $\delta$ sufficiently small.

Hence, we again extend $\Phi$ to $\cl\cC$ by continuity and to the complete phase space $(u,u_x)\in\R^2$ 
by the identity $\Phi|_{(\R^2\setminus\cl\cC)}=\id$. Notice that each $\cC_k$ is isolated and, due 
to this isolation and the integrability of \eqref{eq101} for $f=f_1(u,u_x)$, each $\cC_k$ accepts the 
left/right orientation described in the previous Remark after case {\em (II)}. 

We obtain a Hamiltonian and pendulum realization of \eqref{eq102} for $f=f_1(u,u_x)$ by a 
nonlinearity $g(u)$ defined on the isolated intervals determined by the regions $\cC_k$. 
Then, after the $C^2$ smooth and globally bounded extension of $g(u)$ to $\R$, filling the 
unbounded intervals and the interval gaps, we define the potential $G$ by \eqref{eq313}.
Therefore, due to the partial order \eqref{eq314} restricted to the saddle points on the boundaries of 
$\cC_k$, the phase portraits of \eqref{eq102} for $f=f_0(u,u_x)$ and $f=g(u)$ are qualitatively the 
same. 
Then, preservation of hyperbolicity ensures that $f=f_0(u,u_x)$ and $f=g(u)$ belong to the same 
lap signature class. 
This completes the proof of case {\em (IV)}.

Finally, as a result of \cite{roc12}, Theorem 2, (see also \cite{roc07,firowo12a}), there exists a 
homotopy between the pendulum realization $g(u)$ of $f_0(u,u_x)$ and a pendulum realization 
$\overline g(u)$ of $f_1(u,u_x)$ since, by assumption, both belong to the same lap signature class.
In fact, all global attractors of pendulum realizations with the same lap signature are 
connection equivalent. 

In summary, the global collective homotopy is the composition of: 
\begin{enumerate}[(i)]
\item a Hamiltonian realization homotopy \eqref{eq311}, \eqref{eq312}, and a pendulum realization 
homotopy \eqref{eq303} from $f_0(u,u_x)$ to $g(u)$;  
\item a homotopy between $g(u)$ and $\overline g(u)$;  
\item a reverse pendulum realization homotopy and a reverse Hamiltonian realization homotopy from 
$\overline g(u)$ to $f_1(u,u_x)$.
\end{enumerate}
This concludes the proof of Theorem \ref{th2}. 
$\hfill \square$

\smallskip

\section{Discussion and concluding remarks}

Theorem \ref{th1} is a contribution to the geometric and qualitative theory of dynamical systems 
generated by parabolic partial differential equations. This result, applied to one-dimensional scalar 
semilinear PDEs, extends to the case of reversible nonlinearities the classification of these 
dynamical systems, already initiated in \cite{firowo04,roc07,firowo12a,firowo12b}. 
Similar results are also available in the more complex case of monotone feedback delay differential 
equations. For results and references see \cite{krva11,krva16,nie23}.

\begin{figure}[!ht] \label{fig8}
\begin{center}
\includegraphics[scale = .9]{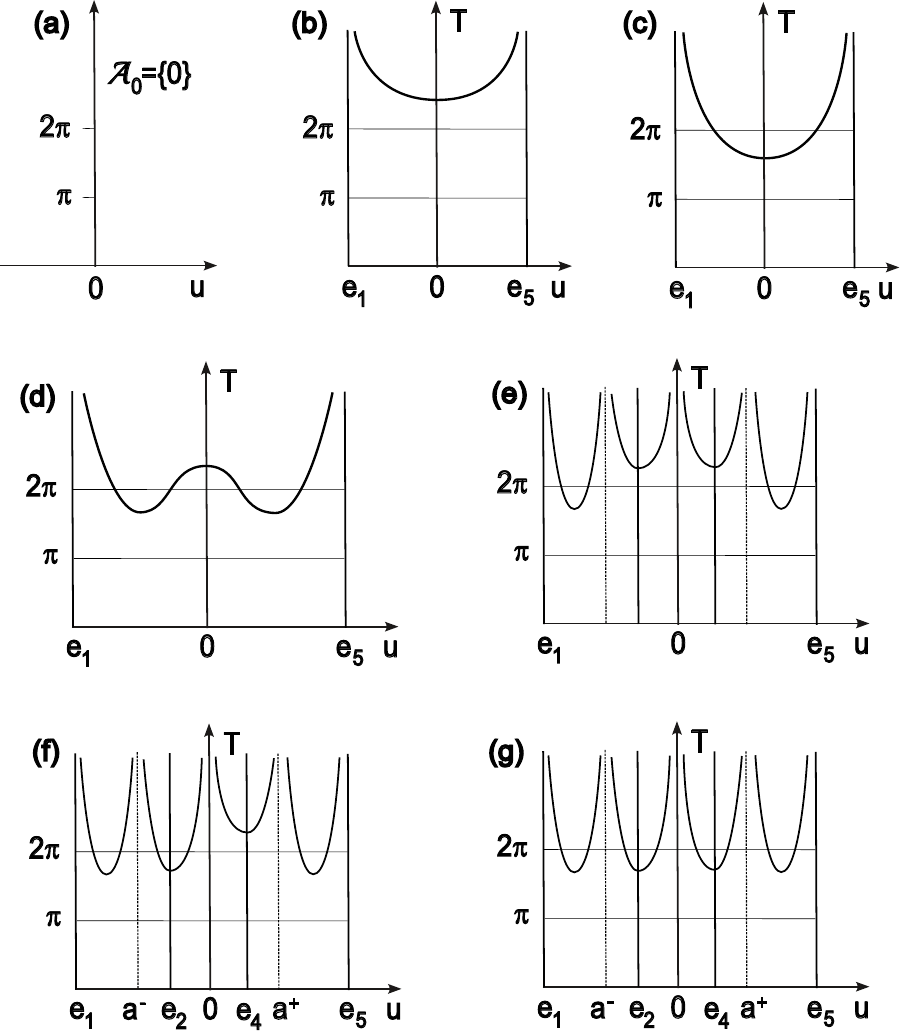}
\end{center}
\caption{\small For the global attractor $\cA_1$ of \eqref{eq101} with lap signature 
$\bigl(\{1,1\}\bigl(\{1\}\bigr)\bigl(\{1\}\bigr)\bigr)$, we consider the family 
$f_s\in\cR, s\in[0,1]$ connecting $\cA_0=\{0\}$ to $\cA_1$. For this example we illustrate the seven 
period maps which satisfy condition (H). The corresponding $s$-intervals are: (a) $0<s<1/7$; 
(b) $1/7<s<2/7$; (c) $2/7<s<3/7$; (d) $3/7<s<4/7$; (e) $4/7<s<5/7$; (f) $5/7<s<6/7$; and, (g) $6/7<s<1$. 
The occurring six bifurcations are discussed in Fig. 9.}
\end{figure} 

Clearly, due to Theorem \ref{th1}, a classification of the global attractors for dynamical systems 
generated by \eqref{eq101} in the reversible class of nonlinearities is provided by the lap signature 
class. As mentioned in Sect. 1, the lap signature class is given by the set of period lap numbers of 
the $2\pi$-periodic orbits endowed with the total order derived from the regular parenthesis structure 
of their nesting in phase space $(u,u_x)$. 
This has the form (see \cite{firowo12a}) 
\begin{equation} \label{eq401}
\Bigl(\{\ell_1^1,\dots,\ell_{k_1}^1\}\Bigl(\{\ell_1^2,\dots,\ell_{k_2}^2\}\Bigl(\dots\Bigr)\dots
\Bigl(\{\ell_1^K,\dots,\ell_{k_K}^K\}\Bigr)\Bigr)\Bigr) \ ,
\end{equation}
where $K$ again denotes the number of cyclicity regions, $k_j$ denotes the number of $2\pi$-periodic 
orbits in the $j^{th}$ annular or punctured disk region (this may be empty, in which case $k_j=0$), and 
$\ell_i^j$ denotes the period lap number of the $i^{th}$ $2\pi$-periodic orbit in the $j^{th}$ cyclicity 
region. The regular parenthesis structure represents the nesting of the periodic orbits. 
To illustrate this representation, we exhibit the lap signature class for the example of Fig. 6 
(with $n=5$, $r=1$, $k_1=2$ and $k_2=k_3=1$) completing the parenthesis order structure shown in 
\eqref{eq315}:
\begin{equation} \label{eq402}
\Bigl(\{1,1\}\Bigl(\{1\}\Bigr)\Bigl(\{1\}\Bigr)\Bigr) \ .
\end{equation}

\bigskip

The present result also shows that in the class of $S^1$-equivariant nonlinearities 
$f^s(u,u_x)=f(s,u,u_x), s\in[0,1]$ 
there is a continuous family of functions $f^s, s\in[0,1]$, connecting the global Sturm attractors 
$\cA_{f^0}=\{{\bf 0}\}$ and $\cA_{f^1}$ through a finite number of local (degenerate) bifurcations. 
Here $s$ denotes the bifurcation parameter. 
Moreover, these bifurcations consist only of pitchfork and Hopf bifurcations. Notice that the 
bifurcations shown are degenerate because the family $f^s\in\cR, s\in[0,1]$, is entirely constructed 
in the symmetry class of reversible nonlinearities $f=f(u,u_x)$, which implies its integrability.
We also illustrate this bifurcation diagram in the example of Fig. 6. We obtain a family 
$f^s, s\in[0,1]$, with six bifurcations points at $s=k/7, 1\le k\le 6$, and seven (degenerate) 
structuraly stable regions $[0,1]\setminus\bigl(\cup_{1\le k\le 6}\{k/7\}\bigr)$. In Fig. 8 we 
illustrate the period maps in these seven regions. Since we are restricted to the integrable case 
$f^s\in\cR$, all the bifurcations are degenerate. For general references see \cite{car81,chha82,guho83}. 

In Fig. 9 we then show the bifurcation diagram from $\cA_0=\{0\}$ to $\cA_1$ and the connection graph 
of a nondegenerate global attractor $\cA_{s>1}$ of \eqref{eq101} homotopic to the (degenerate) $\cA_1$. 
For this we extend the family $f^s$ to $s>1$, connecting the reversible (degenerate) $\cA_1$ to a 
global attractor $\cA_{s>1}$ which is not reversible, see \cite{firowo12a}. In this case, 
$\cA_{s>1}$ is nondegenerate and we have the usual notions of stability of equilibria and 
rotation waves, \cite{firowo12b}. See \cite{firowo04}, Proposition 3.1, and \cite{firowo12b}, 
Theorem 10. 

\begin{figure}[!ht] \label{fig9}
\begin{center}
\includegraphics[scale = 1]{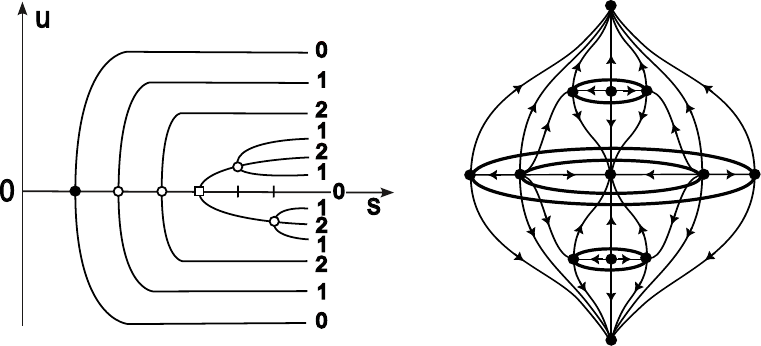}
\end{center}
\caption{\small {\em Left:} The bifurcation diagram for the reversible family $f^s\in\cR, s\in[0,1]$, 
(see Fig. 8). 
The bifurcation points are, respectively: a pitchfork of a stable saddle point (black ball); 
degenerate Hopf bifurcations (white balls); and, a pitchfork of centers (white square). 
Moreover, the sequence of bifurcations leading to $\cA_1$ are: $(i)$ a supercritical pitchfork 
bifurcation of two stable saddle points at $s=1/7$; $(ii)$ a supercritical Hopf bifurcation of an 
unstable $2\pi$-periodic orbit at $s=2/7$; $(iii)$ a supercritical Hopf bifurcation of an unstable 
$2\pi$-periodic orbit at $s=3/7$; $(iv)$ a supercritical pitchfork bifurcation of two centers at 
$s=4/7$: $(v)$ a supercritical Hopf bifurcation of an unstable $2\pi$-periodic orbit at $s=5/7$ in 
the upper branch center; and, $(v)$ a supercritical Hopf bifurcation of an unstable $2\pi$-periodic 
orbit at $s=6/7$ in the lower branch center. At the right margin of the bifurcation diagram we indicate 
the Morse indices corresponding to the equilibria and rotating waves of the nondegenerate Neumann 
section of the global attractor $\cA_{s>1}$. 
{\em Right:} The connection graph of the nondegenerate global attractor $\cA_{s>1}$.}
\end{figure} 

The global attractor $\cA_{s>1}$ has three stable equilibria, denoted by 
$\{{\bf e}_1,{\bf e}_3={\bf 0},{\bf e}_5\}$, and two unstable equilibria, $\{{\bf e}_2,{\bf e}_4\}$, 
where 
\begin{equation} \label{eq403}
\begin{aligned}
\dim W^u({\bf e}_1)&=\dim W^u({\bf e}_3)=\dim W^u({\bf e}_5)=0 \ , \\
\dim W^u({\bf e}_2)&=\dim W^u({\bf e}_4)=3 \ .
\end{aligned}
\end{equation}
In addition, there are four rotating waves, denoted $\{{\bf v}_1,{\bf v}_2,{\bf v}_3,{\bf v}_4\}$ 
(by order of appearance). 
By the normal hyperbolicity of rotating waves, we obtain 
\begin{equation} \label{eq404}
\dim W^u({\bf v}_1)=\dim W^u({\bf v}_3)=\dim W^u({\bf v}_4)=1 \quad , 
\quad \dim W^u({\bf v}_2)=2 \ .
\end{equation}
Therefore, since
\begin{equation} \label{eq405}
\dim \cA_{s>1} = 
\max \Bigl\{\dim W^u({\bf e}_j), 1+\dim W^u({\bf v}_k), 1\le j\le 5, 1\le k\le 4\Bigr\} = 3 \ ,
\end{equation}
the global attractor $\cA_{s>1}$ is three dimensional. The heteroclinic orbit connections follow 
from the connections on the Neumann section of $\cA_{s>1}$, \cite{firowo04,firowo12b}. 
The rotating waves, each one rotating with its own fixed speed $c_j$, appear on the slow stable/unstable 
invariant manifolds of the equilibria $\{{\bf e}_2,{\bf e}_3,{\bf e}_4\}$. 

\bigskip

We conclude by observing that the connection equivalence of global attractors presented in Theorem 
\ref{th1} extends to the general non-reversible case as orbit equivalence, if we assume that all 
spatially nonhomogeneous solutions of \eqref{eq101} are rotating waves, rotating around the circle 
$S^1$ with constant speeds $c_j\ne 0$. 
In this case, we have $\cF=\emptyset$ and all periodic orbits in $\cP$ are rotating waves. So, in 
addition to Theorem \ref{th1}, we can invoke hyperbolicity (H) and the strong Morse-Smale property.

As a final comment, we believe that Theorem \ref{th1} holds for orbit equivalence. To support this 
statement we observe that a small perturbation $h(u,u_x)=\varepsilon v_x$ breaks reversibility and 
the orbit equivalence  statement holds for the homotopy with $f+h$. 

As an application we exhibit the connection graph of the global attractor of \eqref{eq101} for the 
$S^1$-equivariant but non-reversible perturbed Chafee-Infante nonlinearity:
\begin{equation} \label{eq406}
u_t = u_{xx} + \lambda u - u^3 + \varepsilon u_x \quad , \quad x \in S^1 , t \ge 0 \ ,
\end{equation} 
where $\lambda>0$ satisfies $\lambda \ne \lambda_k:=k^2, k\in\N$, and $\varepsilon\ne 0$. 
The set of equilibria of \eqref{eq406}, $\cE=\{{\bf e}_1, {\bf e}_2, {\bf e}_3\}$, has the equilibria 
${\bf e}_1=+\sqrt{\lambda}$ and ${\bf e}_3=-\sqrt{\lambda}$, which are defined and hyperbolic for all 
$\lambda>0$. 
These are saddle points of \eqref{eq102}. The third equilibrium, ${\bf e}_2=0$, is also hyperbolic for 
$\lambda\ne\lambda_1,\lambda_2,\dots$, and is a center of \eqref{eq102}. 

For $\lambda\in(\lambda_k,\lambda_{k+1})$ the set of periodic orbits $\cP=\{v_1,\dots,v_k\}$ has $k$ 
hyperbolic rotating waves which rotate around $S^1$ with speed $\varepsilon$ and $\cF=\emptyset$. 
Hence, the lap signature of the global attractor $\cA_\lambda$ is 
\begin{equation} \label{eq407}
\Bigl( \{1,\dots,k\} \Bigr) \ ,
\end{equation} 
and the connection graph of $\cA_\lambda$ is shown in Fig. 10. As expected, the connection graph of 
this ``spindle attractor'' with multiple periodic orbits in his ``belt'' is a tower, see \cite{leyo21}. 
By the previous results using normal hyperbolicity this connection graph is preserved for 
$\varepsilon=0$, i.e. the reversible case.

\begin{figure}[!ht] \label{fig10}
\begin{center}
\includegraphics[scale = 1]{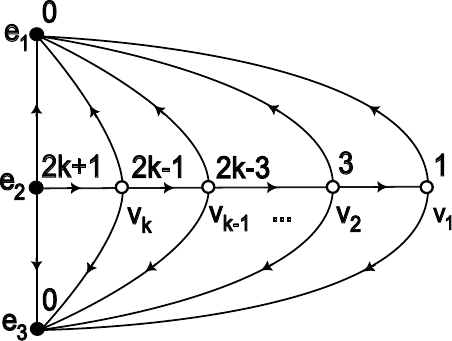}
\end{center}
\caption{\small Connection graph of the global attractor $\cA_\lambda$ of \eqref{eq406} for 
$k^2<\lambda<(k+1)^2$. 
The numbers close to the vertices are the corresponding Morse indices. This shows that 
$\dim\cA_\lambda=2k+1$.}
\end{figure}

As a concluding remark we point out that the nonlinearity $f=f(u,u_x)$ can be regarded as a 
feedback control parameter for the reaction-diffusion equation \eqref{eq101}. This provides a 
very interesting model for certain applications. We mention, for example, the possibility of 
stabilizing certain equilibria corresponding to reaction-diffusion patterns in chemical 
reactors with a proper choice of a nonlinearity $f=f(u,u_x)$.
  
\bigskip \bigskip 

\noindent {\sc Acknowledgments:} 
The author acknowledges the continuous support of Bernold Fiedler (Freie U. Berlin). 
I am deeply thankful for his longtime friendship and collaboration. I am also gratefull 
to the continuous support of Isabel Rocha. Without her constant support, this endeavor 
would not be possible. This work was partially supported by FCT/Portugal 
through the project UIDB/04459/2020 with DOI identifier 10-54499/UIDP/04459/2020.

The author also acknowledges a helpful discussion with Luis Barreira, and the expert comments 
of the three anonymous referees which largely contributed to the improvement of the manuscript.

\end{document}